\documentclass{amsart}

\usepackage{amssymb}
\usepackage{graphicx}
\usepackage{url} 

\newtheorem{theorem}{Theorem}[section]
\newtheorem{lemma}[theorem]{Lemma}
\newtheorem{proposition}[theorem]{Proposition}
\newtheorem{corollary}[theorem]{Corollary}

\newtheorem{definition}[theorem]{Definition}
\newtheorem{example}[theorem]{Example}

\renewcommand{\phi}{\varphi}
\renewcommand{\epsilon}{\varepsilon}
\newcommand{\U}{\mathbb{U}}

\newcommand{\astar}{A^*}

\newcommand{\bdu}{\partial \U}

\newcommand{\bilin}[2]{\langle#1,\,#2\rangle}

\newcommand{\C}{\mathbb{C}}
\newcommand{\calu}{\mathcal{U}}
\newcommand{\calv}{\mathcal{V}}
\newcommand{\calw}{\mathcal{W}}
\newcommand{\clos}{{\rm clos\,}}
\newcommand{\cphi}{C_\phi}
\newcommand{\cpsi}{C_\psi}

\newcommand{\comega}{C_\omega}
\newcommand{\comegah}{C_{\omega,h}}

\newcommand{\ee}[1]{{\rm EE}(#1)}

\newcommand{\goesto}{\rightarrow}

\newcommand{\htwo}{H^2}

\newcommand{\hilb}{\mathcal{H}}

\newcommand{\infsum}{\sum_{n=0}^\infty}

\newcommand{\lhtwo}{\mathcal{L}(\htwo)}

\renewcommand{\notin}{\not\in}

\newcommand{\pointspec}[1]{\sigma_p(#1)}

\newcommand{\phiu}{\phi(\U)}
\newcommand{\psiu}{\psi(\U)}
\newcommand{\phinv}{\phi^{-1}}

\newcommand{\ran}{\mbox{ran\,}}
\newcommand{\re}{\mbox{\rm Re\,}}

\newcommand{\seq}[1]{\{#1\}}
\newcommand{\sstar}{S^*}

\newcommand{\tphi}{T_\phi}

\newcommand{\tomega}{\tilde{\omega}}
\newcommand{\tpsi}{T_\psi}
\newcommand{\tphistar}{\tphi^*}
\newcommand{\tpsistar}{\tpsi^*}

\newcommand{\xstar}{X^*}

\begin{document}






\title[Intertwining Toeplitz operators]{Intertwining relations and extended eigenvalues  for analytic Toeplitz operators}


%

\author{Paul S. Bourdon}
\address{Department of Mathematics, Washington and Lee University, Lexington VA 24450
}

\email{pbourdon@wlu.edu}

\author{Joel H. Shapiro}
\address{Department of Mathematics and Statistics, Portland State University, Portland OR 97207}

\email{shapiroj@pdx.edu}

%

\thanks{}

%

\begin{abstract}
We study the intertwining relation $XT_\varphi=T_\psi X$ where $T_\phi$ and $T_\psi$ are the Toeplitz operators induced on the Hardy space $H^2$ by analytic functions $\varphi$ and $\psi$, bounded  on the open unit disc $\mathbb{U}$, and $X$ is a nonzero bounded linear operator on $H^2$. Our work centers on the connection between intertwining and the image containment  $\psi(\mathbb{U})\subset\phi(\mathbb{U}),$ as well as on the nature of the intertwining operator $X$.  We use our results to study the ``extended eigenvalues'' of analytic Toeplitz operators $T_\varphi$, i.e., the special case $XT_{\lambda\varphi}=T_\varphi X$, where $\lambda$ is a complex number.
\end{abstract}

\maketitle


\section{Introduction}\label{intro} For (bounded) linear operators $A$, $B$, and $X$ on a Hilbert space $\hilb$, the phrase ``$X$ intertwines $A$ and $B$'' means that 
\begin{equation}\label{intertwining}
    XA=BX \quad{\rm and}\quad X\neq 0. 
 \end{equation}  
When it is convenient to deemphasize the intertwining operator $X$, we will simply write (\ref{intertwining}) as ``$A\propto B$''. Warning: the relation $\propto$ is neither symmetric nor  transitive---for a discussion, see Section 2 below.
We will be particularly interested in the special case $\lambda A\propto A$ where $\lambda$ is a complex number, i.e., in the case for which 
 \begin{equation}\label{extended_eigenvalue}
AX=\lambda XA \quad \mbox{(for~some~$X\neq 0$)}.
\end{equation}
 When this happens, $\lambda$ is called an {\em extended eigenvalue\/}  and $X$ an {\em (extended) eigenoperator\/} of $A$.

 The concepts of extended eigenvalue and eigenoperator gained prominence in the 1970s as a result of  Victor Lomonosov's  famous invariant subspace theorem \cite{Lo},  and the following extension of it due independently to Scott Brown \cite{B} and Kim, Moore, and Pearcy \cite{KMP}.

\begin{quote}
{\em If an operator $A$ on a Banach space has a nonzero compact eigenoperator, then $A$ has a nontrivial, hyperinvariant subspace.} 
\end{quote}
The special case where $A$ commutes with a nonzero compact  is the original theorem of Lomonosov.
 
Recently the notion of extended eigenvalue has taken on a life of its own, both in further results of Lomonosov-type \cite{L}, and in studies of the extended eigenvalues and eigenoperators for interesting classes of naturally occurring operators  (\cite{BP, BPL, K, Shk}). In this paper, we continue this latter thread, studying intertwining relationships and extended eigenvalues in the context of function-theoretic operator theory. We work on the most commonly studied Hilbert space of analytic functions, the Hardy space $\htwo$  on the open unit disc $\U$, and consider what is perhaps the most intensely studied class of operators on this space: the {\em analytic Toeplitz operators}.  These are the operators $\tphi:\htwo\goesto\htwo$ defined by pointwise multiplication:
  $$
      \tphi h = \phi h \qquad (h\in\htwo),
  $$
where $\phi$ is a bounded analytic function on $\U$.   

Interest in these operators arises from the special case $\phi(z) \equiv z$,  whereupon $\tphi$ is the forward shift on $\htwo$, whose invariant subspaces are described by the famous theorem of Beurling (see \cite[Section 17.21, page 348]{Rud}, for example). The entire collection of analytic Toeplitz operators is the commutant of this shift (see \cite[Problem 147, page 79]{Hal}, for example).

Our work builds on earlier efforts by James Deddens who, in the early 1970s proved the following 
results, where $\phi$ and $\psi$ are  bounded analytic functions on $\U$.
\begin{itemize}
   \item[{\bf D1}] \cite[Theorem, page 243]{Ded1} $\tphi\propto\tpsi \implies\psi(\U)\subset \clos\phi(\U)$, where ``clos'' denotes ``closure in the complex plane.''
   
   \item[{\bf D2}] \cite[Theorem 2, page 860]{Ded2} {\em If $\phi$ is {\em univalent} then} $$\tphi\propto\tpsi \iff \psi(\U)\subset \phi(\U).$$
\end{itemize}
According to Theorem D1, if $\tphi\propto\tpsi$ then each value of $\psi$ lies either in $\phi(\U)$ or on its boundary.  We show here that more is true: {\em the values of $\psi$ that lie on the boundary of $\phi(\U)$ must be ``capacitarily isolated''} in the sense that some neighborhood of such a value must intersect $\partial\phi(\U)$ in a set of logarithmic capacity zero. If, in particular, the boundary of $\phi(\U)$ consists entirely of nondegenerate continua, then we must have $\psi(\U)\subset\phi(\U)$. This recovers, as a special case,  the forward implication of Theorem D2. We generalize the reverse implication of Theorem D2 in a different way:  If $\phi$ is a {\em covering map\/} and $\psi(\U)\subset\phi(\U)$, then  $\tphi\propto\tpsi$. So, for example, our generalizations combine to yield:
\begin{quote} {\em
If  $\phi$ is a covering map taking $\U$ onto a domain whose boundary consists entirely of nondegenerate continua, then $$\tphi\propto\tpsi \iff \psi(\U)\subset\phi(\U).$$
}
\end{quote}
Our results on capacitary isolation and covering maps occupy Section 3, with the necessary background material and preliminary results collected in Section
2.

 It is known that the containment $\psiu\subset\phiu$ is not, in general, sufficient to guarantee that $\tphi\propto\tpsi$.  As a corollary of a theorem connecting the equality of the commutants of $A$ and $A^n$ to the unimodular extended eigenvalues of $A$,   Carl Cowen \cite[p.\ 339]{Cow2} shows that if  $f$ is a Riemann map of $\U$ onto the slit disc $\U\setminus (-1, 0]$ and $\phi = f^2$ while $\psi = -f^2$, then $\tphi \not\propto \tpsi$ even though $\psi(\U) = \phi(\U)$.   Cowen's example settles one-half of a conjecture due to Deddens  \cite[page 244]{Ded1}, who speculated that the condition ``\,$\overline{\psi(\U)}\subset$ point spectrum of $\tphi^*$\,''  is necessary and sufficient for $\tphi\propto\tpsi$ (here the ``overline'' means ``complex conjugate'').    Since the point spectrum of $\tphi^*$ always contains $\overline{\phi(\U)}$, Cowen's example settles in the negative the sufficiency part of Deddens's conjecture.  Necessity remains open.  
 
    Most intertwining issues are trivial when either $\phi$ or $\psi$ is a constant function; for example, the necessity part of Deddens' conjecture is easily seen to be valid.  Suppose that $\psi$ is the  constant function $z\mapsto \lambda$ and that  $\tphi$ and $\tpsi$ are intertwined by $X$; then $\tphistar X^* = \bar{\lambda}X^*$, which, since $X\ne 0$, tells us $\bar{\lambda}$ is a eigenvalue of $\tphistar$. On the other hand,   if $\lambda$ is the constant value of $\phi$ and $\tphi\propto \tpsi$, then $\psi =\phi$ so that both operators act as $\lambda$ times the identity on $\htwo$, and, once again, the necessity part of Deddens' conjecture holds. 

 In Section 4, we take up the issue of characterizing intertwining operators for analytic Toeplitz operators, establishing that, e.g.,  if $X\tphi = \tpsi X$ and $\phi$ is univalent, then $X$ must be a weighted composition operator. More generally, we show that if $\phi$ singly covers some nonvoid open set $W\subset \phi(\U)$ and $\psi(\U)\cap W \ne \emptyset$, then $X$ must be a weighted composition operator.  This result leads to further examples illustrating that the sufficiency part of Deddens' conjecture is invalid.  For instance, we show that if for some $p>1$, the function $\phi$ is  $p$-valent and $\phi$ singly covers some nonvoid open subset of $\phi(\U)$ while $\phi-\phi(z_0)$ has a zero of order $p$ for some $z_0\in\U$, then there exist bounded analytic functions $\psi$ with $\psiu\subset \phiu$ for which $\tphi \not\propto \tpsi$.   Our work provides a new proof  that when $\phi$ singly covers a nonvoid open subset of its image, then $\{\tphi\}' = \{T_z\}'$. The original proof is due to Deddens and Wong \cite[Corollary 4]{DW}. 
 
 A linear combination of operators that intertwine $A$ and $B$ will also intertwine $A$ and $B$.  We use this observation in Section 4 to show that if $\phi$ is not univalent, there will exist a $\psi$ and an $X$ intertwining $\tphi$ and $\tpsi$ such that $X$ is not a weighted composition operator  but rather a linear combination of weighted composition operators.  On the other hand, we show that if $\tphi$ and $\tpsi$ are intertwined by a linear combination of weighted composition operators, then they are intertwined by a single composition operator.  Thus, the collection of pairs of analytic Toeplitz operators that may be intertwined by composition operators is not enlarged by allowing intertwining by linear combinations of weighted composition operators.
 
 In Section 5, we consider the relation $\tphi\propto \tpsi$ when $\phi$ is an inner function.   In this context, we provide examples showing that sometimes intertwining operators cannot be composition operators, (and hence cannot be linear combinations of weighted composition operators).  
We rely on, and provide a new proof for, the following result due to  Deddens \cite{Ded2}:   $\tphi\propto\tpsi$ whenever $\phi$ is inner and $\psi$ is a self-map of $\U$.

 In Section 6,  we consider the relation  $\tphi\propto\tpsi$ when $\psi$ is an inner function.  We show that if $\tphi\propto\tpsi$ and $\psi$ is a nonconstant inner function, then the point spectrum of $\tphistar$ must contain $\U$.   Thus, the necessity part of Deddens' conjecture holds when $\psi$ is inner.   Also,  this result shows that $\tphi\not\propto \tpsi$ whenever $\psi$ is inner and $\phi$ is outer.  (It is easy to produce examples where $\tphi \propto \tpsi$ where $\phi$ is inner and $\psi$ is outer.)     If $\psi$ is inner, then $\tpsi$ is a completely nonunitary isometry, i.e.,  a unilateral shift.  Our work in Section 6 is carried out at the more general level of isometries: we prove that if $S$ is a completely nonunitary isometry on a Hilbert space $\hilb$ and $A:\hilb\rightarrow \hilb$ is a bounded linear operator, such that $XA = SX$ for some nonzero operator $X$ on $\hilb$, then every point of $\U$ is an eigenvalue of $A^*$.  
  
  In the final section of this paper, we consider  consequences of our results and those of Deddens'  to the study of extended eigenvalues of analytic Toeplitz operators.   When $\phi$ is a covering map, the geometry of $\phi(\U)$ quickly determines the collection $\ee{T_\phi}$ of extended eigenvalues of $T_\phi$.   On the other hand, our work in Section 4 provides examples showing that, in general, the determination of $\ee{T_\phi}$ is a challenging problem.

\section{Preliminaries}\label{preliminaries} 

Crucial to our work will be the class of {\em weighted composition operators}  on $\htwo$. If $\omega$ is holomorphic on $\U$ with $\omega(\U)\subset\U$, then the  composition operator $\comega$ defined by
$$
     \comega f = f\circ\omega \qquad (f\in\htwo)
$$
is, by Littlewood's Theorem (see \cite[Chapter 1]{Sh}, for example), a bounded linear operator on $\htwo$.  If $\omega$ is a holomorphic self-map of $\U$ and $h\in H^2$ is such that the product $h\cdot (f\circ \omega)$ belongs to $\htwo$ whenever $f\in \htwo$, then the Closed Graph Theorem shows the weighted composition operator $\comegah$ defined by $\comegah f = h\cdot (f\circ\omega)$  is a bounded linear operator on $\htwo$.  Observe that $C_{\omega, h}$ is clearly bounded on $H^2$ whenever $h$ is bounded on $\U$.  However, the boundedness of $h$ isn't necessary; for example if $h(z)= \log(1-z)$ and $\omega(z) = z/2$ then $\comegah$ maps $H^2$ into itself.  
 
 The connection with analytic Toeplitz operators is this: if $\phi$ is, as usual, a bounded analytic function on $\U$, then $\comegah$ intertwines $\tphi$ and $T_{\phi\circ \omega}$; that is 
\begin{equation} \label{composition_intertwining}
     \comegah  \tphi = T_{\phi\circ\omega}\comegah. 
\end{equation}
Here, we are assuming $h$ is not the zero function so that $\comegah$ is not the zero operator.

\begin{example}  $T_{2z}\propto T_z$ with $X = C_{z/2, 1} = C_{z/2}$  serving as an intertwining map.  On the other hand, Deddens' Theorem D1 shows that $T_z\not\propto T_{2z}$.  Thus, the relation $\propto$ is not symmetric. 
\end{example}

\noindent
{\em Remarks.} (a) It is easy to give a simple, direct proof that $T_z\not\propto T_{2z}$.  Suppose $XT_z = T_{2z}X$ for some nonzero operator $X$ on $\htwo$.   Choose $g\in H^2$ with $\|g\| =1$  such that $\|Xg\| > \|X\|/2$.  Then 
$$
\|X\| < 2\|Xg\| = \|T_{2z}X g\|  = \|XT_z g\| \le \|X\| ,
$$
a contradiction.  

(b)   In the finite-dimensional setting, however, the relation $\propto$ {\em  is\/} symmetric.   
\begin{proposition}  Suppose that  $A$ and $B$ are $n\times n$ matrices.  If there is a nonzero $n\times n$ matrix $X$ such that $XA = BX$, then  there is a nonzero $n\times n$ matrix $Y$ such that $YB = AY$.  
\end{proposition}
\begin{proof}    We apply the Sylvester--Rosenblum Theorem (see, e.g., \cite{BR}), which asserts that if $A$ and $B$ are operators on a Banach space with disjoint spectra, then the equation $AX-XB=Y$ has a unique solution for every operator $Y$.  Suppose that  there is a nonzero $n\times n$ matrix  such that $XA = BX$.  By the Sylvester--Rosenblum Theorem, there is  a number $\lambda\in \sigma(A)\cap \sigma(B)$.  Let $a$ be an eigenvector for $A$ with corresponding eigenvalue $\lambda$.   Let $b$ be an eigenvector for $B^T$ with corresponding eigenvalue $\lambda$.   Let $Y=ab^T$, a nonzero $n\times n$ matrix. Then 
    $$
      AY = Aab^T = \lambda ab^T =\lambda Y
    $$  
    and
    $$
      B^TY^T = B^Tba^T = \lambda ba^T =\lambda Y^T
    $$
    so that $YB = \lambda Y = AY$, as desired. 
    \end{proof}

  The proof of the preceding proposition shows that if the $n\times n$ matrices $A$ and $B$ have a common eigenvalue, then $A\propto  B$.   This observation, along with the Sylvester--Rosenblum Theorem provides examples illustrating that the relation $\propto$ is not transitive.   Let $A$, $B$, and $C$ be $n\times n$ matrices such that $\sigma(A)\cap \sigma(B)\ne \emptyset$, $\sigma(B)\cap\sigma(C)\ne \emptyset$, and $\sigma(A)\cap\sigma(C)= \emptyset$.  The $A\propto B$ and $B\propto C$, yet $A\not\propto C$.

We now present several known results  needed in the sequel. For the sake of completeness, we sketch out some of the proofs.

An important class of test functions for analytic Toeplitz operators is the set of {\em reproducing kernels.} For each point $a\in\U$, the ($\htwo$) reproducing kernel for $a$ is the function 
  $$
      K_a(z) = \frac{1}{1-\overline{a}z} \qquad (z\in\U).
  $$
  The importance of these functions for our study comes from the following well-known result.
\begin{proposition}\label{tphistar_pointspectrum}
$\tphi^* K_a=\overline{\phi(a)}K_a$    for each $a\in\U$. 
\end{proposition}
\begin{proof}
  For each $f\in\htwo$,
  $$
    \bilin{f}{\tphistar K_a} = \bilin{\tphi f}{K_a} = \bilin{\phi f}{K_a} 
                                    = \phi(a)f(a) = \bilin{f}{\overline{\phi(a)}K_a}.
  $$
\end{proof}

\noindent
Note that the preceding proposition shows $\|\tphi\| \ge \|\phi\|_\infty : = \sup\{|\phi(z)|: z\in \U\}$.  On the other hand, the integral representation of the norm of $\htwo$ makes it clear that $\|\tphi\| \le \|\phi\|_\infty$;  thus, $\|\tphi\| = \|\phi\|_\infty$.

\begin{proposition} \label{first_containment}
{\rm\cite[page 273]{Ded1}}\ \  $\tphi\propto\tpsi \implies \psiu\subset\clos{\phiu}$
\end{proposition}

\begin{proof}
The hypothesis is equivalent to saying that there exists a bounded operator $X\neq 0$ on $\htwo$ such that $\tphistar X^* = X^*\tpsistar$. Upon applying both sides of this equality to reproducing kernels $K_a$, we find from Proposition \ref{tphistar_pointspectrum} that either: (i)  $X^*K_a=0$, which can happen at worst on a discrete sequence of points in $\U$ (in fact, at worst on a Blaschke sequence), or (ii) $\overline{\psi(a)}$ is an eigenvalue of $\tphistar$, and so belongs to the spectrum of $\tphistar$, which, by Wintner's Theorem (see \cite[Problem \#247, page 139]{Hal}, for example) is the closure of $\overline{\phiu}$.
Thus, if $a\in\U$ then there is a sequence $a_n$ of points in $\U$ converging to $a$ with $\psi(a_n)$ in the closure of $\phiu$ for each $n$. Thus, $\psi(a)$ is in the closure of $\phiu$. 
\end{proof}

\begin{proposition}\label{pointspec}
$\overline{\alpha}\in\pointspec{\tphistar}$ $\iff$ $\phi-\alpha$ has a nontrivial inner factor.
\end{proposition}

\begin{proof}
$\overline{\alpha}\in\pointspec{\tphistar}$ means that the null space of $T_{\phi-\alpha}^*$ is nontrivial, i.e., that the range of $T_{\phi-\alpha}$, which equals $(\phi-\alpha)\htwo$, is not dense in $\htwo$. By Beurling's Theorem \cite[Theorem 17.21, page 384]{RF}, this is equivalent to saying that $\phi-\alpha$ has no nontrivial inner factor.
\end{proof}

\begin{proposition}\label{strongconv}  Suppose that $\phi$ is an analytic self-map of $\U$.  Then $(\tphistar)^n$ converges to $0$ strongly.
\end{proposition}
\begin{proof} For any $a\in \U$,   
$(\tphistar)^n K_{a} = \overline{\phi(a)}^{\,n} K_{a} \rightarrow 0$ as $n\rightarrow \infty$ because $|\phi(a)| <1$.   Thus, if $f$ is a linear combination of reproducing kernels, $\lim_n \|(\tphistar)^n f\| = 0$.  The proposition now follows  because the linear span of  $\{K_\alpha: \alpha\in \U\}$ is dense in $\htwo$ and  $\|\tphistar\| = \|\tphi\| \le \|\phi\|_\infty \le 1$.\end{proof}

  Let $N^+$ denote the  {\it Smirnov\/} class \cite[Section 2.5]{Dur}, and recall that $H^2\subset N^+$.
   \begin{proposition} \label{rfthm}
{\rm (The ``Rudin--Frostman Theorem'' \cite{RF}).} Suppose that $f\in N^+$. Then $f-\alpha$ has no singular inner factor for all $\alpha\in\C$ outside of a (possibly empty) set of logarithmic capacity zero.
\end{proposition}

\section{Covering and Isolation}\label{equivalent}
The work of this section shows that intertwining is equivalent to image-containment whenever $\phi$ is a covering map whose image has a ``nice'' boundary. We begin with an example showing that some restriction on the boundary is needed.
\begin{example}\label{sing_inner_example} {\rm According to Proposition \ref{first_containment}, if $\tphi\propto\tpsi$ then for each $a\in\U$ the point $\psi(a)$ lies either in $\phiu$ or on $\partial\phiu$.
To see that the latter possibility can actually happen, consider the pair of functions $\psi(z)\equiv z$ and $\phi$ the ``unit singular function''
   \begin{equation}  \label{unit_singular}
    \phi(z)= \exp\left\{\frac{z+1}{z-1}\right\}\qquad (z\in\U),
    \end{equation}
an inner function with $\phiu=\U\setminus \{0\}$. Then $\psi(0)=0\not\in\phiu$, 
but according to a result of Deddens \cite[Cor. 2, page 861]{Ded2}, which we will discuss in some detail in \S\ref{inner_tphi} (see Theorem \ref{intertwining_inner}), $\tphi\propto\tpsi$. Note that in this example the point $\psi(0)=0$, while not belonging to $\phiu$, is {\em isolated\/} in $\partial\phiu$. }
\end{example}

The following result shows that in situations like this some kind of isolation is inevitable.
 To state it concisely let's agree to call a point $w$ of a plane set $S$ {\em capacitarily isolated in $S$\/} if  some neighborhood of $w$ intersects $S$ in a set of logarithmic capacity zero.
\begin{theorem}\label{capisothm}
 If $\tphi \propto \tpsi$ then, for every $a\in\U$, $\psi(a)$ is either in $\phiu$ or is capacitarily isolated in $\partial\phiu$. 
\end{theorem}
\begin{proof}
Fix $a\in\U$ and suppose $\psi(a)$ is not in $\phiu$, and so, by Proposition \ref{first_containment},  lies on $\partial\phiu$. Suppose for the sake of contradiction that $\psi(a)$ is {\em not\/} capacitarily isolated in $\partial\phiu$. Then for each positive integer $n$ the open disc of radius $1/n$ centered at $\psi(a)$ intersects $\partial\phiu$ in a set of positive capacity, so by the Rudin--Frostman Theorem (Proposition \ref{rfthm}) this intersection contains a point $w_n$ such that $\phi-w_n$ has no singular inner factor. Since $w_n\notin\phiu$ this means that $\phi-w_n$ is {\em outer\/}, so by Proposition \ref{pointspec}, $\overline{w_n}\notin \pointspec{\tphi^*}$.

Now, $\psiu$ is an open set and $w_n\goesto \psi(a)$, so we may without loss of generality assume that {\em every\/} point $w_n$ lies in $\psiu$. By elementary function theory, we can choose a sequence $(z_n)$ of points in $\U$ with $z_n\goesto a$ and, for each $n$,  $\psi(z_n)=w_n$.  

Finally, since $\tphi\propto\tpsi$ there exists a nonzero operator $X\in\lhtwo$ for which $X\tphi=\tpsi X$, i.e., 
$$
\tphistar\xstar=\xstar\tpsistar \qquad{\rm with}~\xstar\neq 0.
$$ 
Thus, for every $n$ we have
from Proposition \ref{tphistar_pointspectrum}:
\begin{equation} \label{reproker_intertwining}
    \tphi^*(X^*K_{z_n}) = X^*\tpsi^*K_{z_n} =\overline{\psi(z_n)}(X^*K_{z_n}),
\end{equation}
so if $X^*K_{z_n}$ were $\neq 0,$ then $\overline{\psi(z_n)}$ would be a $\tphi^*$-eigenvalue. But, for each $n$, $\psi(z_n)=w_n$ was chosen so that  it is {\em not\/} a $\tphi^*$ eigenvalue. Thus, $X^*K_{z_n}=0$ for each $n$. Since the sequence $(z_n)$ converges to $a\in\U$, the linear span of $\{K_{z_n}\}_1^\infty$ is dense in $H^2$, and thus we must have $X^*=0$, which contradicts the fact that $X\neq 0$. Thus, $\psi(a)$ is capacitarily isolated in $\partial\phiu$.
\end{proof}

Because any nondegenerate  continuum has positive logarithmic capacity, the following is a consequence of Theorem~\ref{capisothm}.  
\begin{corollary} \label{continua_corollary}
If $\partial\phiu$ consists entirely of nondegenerate continua, then
  $$
      \tphi\propto\tpsi \implies \psiu\subset\phiu.
  $$
\end{corollary} 

Example \ref{sing_inner_example} shows that some restriction on $\partial\phiu$ is needed in order for intertwining to imply image containment. In the corollary above, we could as well have made the weaker assumption that $\partial\phiu$ is ``capacitarily perfect'' in the sense that no point is capacitarily isolated.

Corollary \ref{continua_corollary} has a converse for covering maps, with no extra hypotheses needed on $\partial\phiu$.

\begin{theorem}\label{covering_map}
   If  $\phi$ is a covering map, then 
   $$
      \psiu\subset\phiu \implies \tphi\propto\tpsi.
   $$ 
 \end{theorem}
\begin{proof}  
As we observed in the Introduction, Deddens \cite{Ded2} proved this result for the special case of $\phi$ univalent. For this case, he defined a holomorphic function $\omega:\U\goesto\U$ by $\omega = \phi^{-1}\circ\psi$, so that 
\begin{equation} \label{subord}
\psi = \phi\circ\omega, 
\end{equation}
hence $C_\omega\tphi  = \tpsi C_\omega$ by (\ref{composition_intertwining}),
which completes the proof, with the intertwining operator $X=C_\omega \neq 0$.

For $\phi$ a covering map,  there still exists a holomorphic ``lift'' $\omega:\U\goesto\U$ satisfying (\ref{subord})  (see \cite[Corollary 2.6, page 114]{Con}), so the proof goes through as above. 
\end{proof}

In the covering-map case, the idea behind the existence of the lift $\omega$ of $\psi$  is roughly this: Since $\phi$ is a local homeomorphism, $\omega$ exists locally, and since $\phi$ is in fact a covering map, any local element of $\omega$ can be continued arbitrarily in $\U$, hence, by the Monodromy Theorem, to all of $\U$.

Putting together Corollary \ref{continua_corollary} and Theorem \ref{covering_map}, we obtain the following corollary.

\begin{corollary} \label{nscontainment}
Suppose $\phi$ is a covering map with $\partial\phiu$ a union of nondegenerate continua. Then 
  $$
       \tphi \propto \tpsi \iff \psiu\subset\phiu.
  $$
\end{corollary}

\section{Intertwining Operators} \label{counterexsec}
  We continue to assume that $\phi$ and $\psi$ are  bounded analytic functions on $\U$.   
  \begin{definition}
  We say that $\psi$ is {\em subordinate to} $\phi$ provided there is a holomorphic self-map $\omega$ of $\U$ such that $\psi = \phi\circ \omega$. 
  \end{definition}
We caution the reader that this definition of subordinate differs from the classical one in that we do not require $\omega(0) = 0$.   Our slightly less restrictive notion of subordination has already appeared at the beginning of Section \ref{preliminaries}, and in the proof of Theorem \ref{covering_map}, the point being that if $\psi=\phi\circ\omega$ then $\cphi$ and $\cpsi$ are intertwined by $\comega$.

  In this section, we are interested in the converse implication---establishing hypotheses on $\phi$ and $\psi$ such that $\tphi\propto\tpsi$ implies $\psi$ is subordinate to $\phi$.     Our work will show that only weighted composition operators can intertwine $\tphi$ and $\tpsi$ when $\phi$ is univalent.  It will also provide additional examples illustrating that the containment $\psiu\subset\phiu$ is not sufficient to ensure $\tphi\propto\tpsi$.  
  
We say that ``$\phi$ singly covers a set $W$'' whenever $\phi$ takes some set $V\subset \U$ univalently onto $W$ with $\phinv(W)=V$, i.e., the points of $V$ are the {\em only ones\/} in $\U$ taken by $\phi$ into $W$.

\begin{theorem}\label{continuation_theorem}
Suppose $\phi$  singly covers a nonvoid open subset  of $\phiu$  that has nontrivial intersection with $\psiu$.  Then
     \begin{center}
          $\tphi \propto \tpsi ~~\implies ~~\psi$ is subordinate to $\phi$,          
     \end{center}
 in particular, $\psi(\U)\subset\phi(\U)$.
\end{theorem}

  \begin{proof}{\it of Theorem~\ref{continuation_theorem}} \quad  Suppose that $\psi$ is constant, taking the value $c$ at each $z\in\U$; then $c$ must belong to the region singly covered by $\phi$, and $\psi$ is subordinate to $\phi$: $\psi = \phi\circ \omega$,  where $\omega$ is a constant function whose  value is the preimage of $c$ under $\phi$. (Thus, when $\psi$ is constant, the additional hypothesis that $\tphi \propto \tpsi$ is not needed to obtain the conclusion of the theorem.)  For the remainder of the proof, we assume that $\psi$ is nonconstant.  
  
  Suppose that $\tphi \propto \tpsi$ so that there is a nonzero operator $X$ on $\htwo$ that intertwines $\tphi$ with $\tpsi$.  Consider the map $F:\U\goesto\htwo$ defined by
   $$
       F(z) = \xstar K_z \qquad (z\in\U),
   $$
   where $K_z$ is the $\htwo$-reproducing kernel for the point $z$. Then $F$ is a holomorphic $\htwo$-valued map on  $\U$, which, thanks to the fact that $\xstar\neq 0$, is not identically zero. Thus, the set
   $$
       Z(F) = \seq{z\in\U: F(z)=0}
   $$
is at worst a sequence in $\U$ that tends to $\bdu$ (in fact, it is at worst a Blaschke sequence).  As in the proof of Theorem \ref{capisothm}, we see from Proposition \ref{tphistar_pointspectrum} that $T_\phi^*F(z) = \overline{\psi(z)}F(z)$, so whenever  $z\notin Z(F)$ the vector $F(z)$ is a $\tphistar$-eigenvector for the eigenvalue $\overline{\psi(z)}$. 

Referring back to the hypotheses of our theorem, let  $W$ be the nonvoid open subset of $\phi(\U)$ that is singly covered by $\phi$ and let $V= \phi^{-1}(W)$. Thus, $\phi$ maps $V$ univalently onto $W$ and  $\phi^{-1}$ is a  holomorphic map from $W$ onto $V$.  Let $G=\psi^{-1}(\psi(\U)\cap W)$; then $G$ is open, nonvoid (since $\psi(\U)\cap W\neq\emptyset$), and $\psi(G)\subset W$.  

\bigskip\noindent
{\sc Claim:} For each $z\in G\setminus Z(F),$ the $\tphistar$-eigenvalue $\overline{\psi(z)}$ has {\em multiplicity one.} 

\bigskip\noindent
Granting this for the moment, note that because $\psi(G) \subset W$, the function $\omega := \phinv\circ\psi:G\goesto V$ is holomorphic. Observe that, by Proposition \ref{tphistar_pointspectrum}, for each $z\in G$,
  $$
     \tphistar K_{\omega(z)} = \overline{\phi(\omega(z))} K_{\omega(z)}
                                        = \overline{\psi(z)} K_{\omega(z)}.
  $$
Thus, for each $z\in G\setminus  Z(F)$ both $K_{\omega(z)}$ and $F(z)$ are $\tphistar$-eigenvectors for the eigenvalue $\overline{\psi(z)}$, which, by the {\sc Claim}, has multiplicity one.  Thus, for each such $z,$ there exists a nonzero complex number, which we will call $\overline{h(z)}$, such that $F(z)=\overline{h(z)}K_{\omega(z)}$. Let's set $h\equiv 0$ on $Z(F)$, so that
\begin{equation} \label{h_via_reproker}
  F(z) = \overline{h(z)}K_{\omega(z)} \qquad (z\in G).
\end{equation}
 
 To finish the proof (modulo proving the {\sc Claim}) we examine the nature of the function $h$ and the operator $X$. Observe that for each $f\in\htwo$ we have $Xf\in\htwo$, and for each $z\in G$:
   \begin{eqnarray*}
      Xf(z) &=& \bilin{Xf}{K_z} =\bilin{f}{\xstar K_z}  \\
              &=& \bilin{f}{F(z)} = \bilin{f}{\overline{h(z)}K_{\omega(z)}}\\
              &=& h(z) \bilin{f}{K_{\omega(z)}} = h(z) f(\omega(z)), 
   \end{eqnarray*}
where (\ref{h_via_reproker}) is used in the second equality of the second line.  Thus,
 \begin{equation} \label{extension_equation} 
       Xf=h\cdot (f\circ\omega) \quad {\rm on}~G.
 \end{equation}
 Since $Xf$ is holomorphic on $\U$, equation (\ref{extension_equation}) shows that product  $h\cdot (f\circ\omega)$ is not only holomorphic on $G$, it has a holomorphic extension to $\U$. Upon setting $f\equiv 1$ on $\U$ we see from (\ref{extension_equation}) that $h$ itself extends holomorphically to $\U$, while the assignment $f(z)\equiv z$ shows that $h\cdot\omega$ does the same. Thus,         $\omega$ extends to a quotient of functions holomorphic on $\U$, i.e., to a function meromorphic on $\U$.  Call this extension $\tomega$, and note that its only possible singularities are zeros of $h$ in $\U\setminus G$.  

We claim that $\tomega$ extends holomorphically to $\U$, and that this extension maps $\U$ into itself.  
To see this,
let $S = \{z\in \U: h(z) = 0\}$,  let $n$ be an arbitrary positive integer, and let $f_n$ be defined on $\U$ by $f_n(z) = z^n$.  Take $f= f_n$ in (\ref{extension_equation})  to see that $\omega^n = (Xf_n)/h$ on $G\setminus S$.     For each $z\in G\setminus S$, we have
$$
\tomega(z)^n = \omega(z)^n = \frac{(Xf_n)(z)}{h(z)}
$$
and it follows that the functions $\tomega^n$ and $(Xf_n)/h$, which are analytic on $\U\setminus S$, must agree on $\U\setminus S$.
We therefore have for every point $a\in \U\setminus S$ and every positive integer $n$, 
\begin{equation}\label{kineq}
|h(a)||\tomega(a)|^n  = |(Xf_n)(a)| = |\langle Xf_n, K_a\rangle| \le \|Xf_n\|\|K_a\| \le \frac{\|X\|}{\sqrt{1-|a|^2}},
\end{equation}
where the final inequality holds because, for each positive integer $n$, the function $f_n$ has unit norm in $\htwo$.  
Now, take $n$-th roots of the leftmost and rightmost quantities in (\ref{kineq}), and take the limit as $n\rightarrow\infty$ to obtain (because $h(a)\neq 0$) $|\tomega(a)| \le 1$.   Hence, $\tomega$ is bounded on $\U\setminus S$ so that its singularities are removable. 

The holomorphic extension of $\tomega$ to $\U$, which we still label $\tomega$,  satisfies $|\tomega(z)| \le 1$ for all $z\in \U$.  Because $\tomega$ is nonconstant, the open-mapping theorem shows that actually $|\tomega(z)| < 1$ for $z\in U$; that is, $\tomega$ is a holomorphic self-map of $\U$, as desired.
Thus, $\phi\circ\tomega$ is a function holomorphic on $\U$ that agrees with $\psi$ on the nonvoid open set $G$, and therefore agrees with $\psi$ everywhere on $\U$.  In other words, $\psi$ is subordinate to $\phi$, as we wished to show.

 It remains to prove the {\sc Claim}. This will follow from our showing that for $b\in W$
   $$
     \ker(\tphistar -\overline{b}I)=\ker T_{\phi-b}^* = (\ran T_{\phi-b})^\perp
  $$
has dimension one, i.e., that the closure of $\ran(T_{\phi-b})$ has codimension one.  Our hypotheses on $\phi$ guarantee that $\phi-b$ takes the value zero at exactly one point $a\in\U$, that $\phi'(a)\ne 0$, that $a$ is necessarily in $V$, and that if $r>0$ is chosen so that $V$ contains the open disc $D$ centered at $a$ with radius  $r$, then $\phi^{-1}(\phi(D)) = D$.  Replacing $r$ with $r/2$ if necessary, we may assume that $D$ is contained in a proper subdisc of $\U$.  We have
  $$
      \phi(z) - b = (z-a)g(z)
  $$
for $z\in\U$, where $g$ is a nonvanishing holomorphic function on $\U$ that is necessarily bounded on $\U$ (since $\phi$ is).   Furthermore, $g$ is bounded away from $0$ near $\partial\U$; otherwise  there is a sequence $(z_n)$ in $\U$ approaching $\partial \U$ such that $(\phi(z_n))$ approaches $b$, and one could choose $m$ large enough so that both $z_m\not\in D$ and $\phi(z_m)\in \phi(D)$, a contradiction.  Thus, $T_g$ is invertible and we have
  $$
      \ran{T_{\phi-b}} =(\phi-b)\htwo = (z-a)g\htwo = (z-a)\htwo,
  $$
Thus, the range of $T_{\phi-b}$ is already closed, and it has codimension one, which validates our claim and completes the proof of the theorem.
\end{proof}

Suppose that $\phi$ and $\psi$ satisfy the hypotheses of the preceding theorem and that $X$ intertwines $\tphi$ and $\tpsi$.  By the proof of the theorem, the function $\omega := \phi^{-1}\circ\psi$, initially defined on a nonvoid open subset $G$ of $\U$, extends to a holomorphic self-map $\tomega$ of $\U$.  In addition,  there is a function $h\in \htwo$, such that for every $f\in\htwo$,  $Xf = h\cdot(f\circ\omega)$ on $G$. Observe that $C_{\tomega, h} f$ agrees with $Xf$ on $G$ and hence on all of $\U$.   Thus, $X$ is the weighted composition operator $C_{\tomega, h}$.  

  If $\phi$ is univalent and $\tphi\propto \tpsi$, then by Deddens' Theorem D2 (discussed in \S 1) $\psiu\subset\phiu$, and $\phi$ and $\psi$ clearly satisfy the hypotheses of Theorem~\ref{continuation_theorem}.  Thus, by the discussion of the preceding paragraph, we have the following corollary. 
\begin{corollary} If $\phi$ is univalent  and $X$ intertwines $\tphi$  and  $\tpsi$,  then $X$ is a weighted composition operator.
\end{corollary}

In general, intertwining maps for analytic Toeplitz operators need not be weighted composition operators. In fact:
\begin{quote}
{\em As long as $\phi$ is not univalent, there will exist a $\psi$ and an $X$ such that $X$ intertwines $\tphi$ and  $\tpsi$ and $X$ is not a weighted composition operator.}
\end{quote}

The  argument is simple.   First, if $\phi$ is constant and $\tphi\propto \tpsi$, then as we pointed out in the Introduction, $\tphi$ and $\tpsi$ equal $\lambda I$, where $\lambda$ is the constant value of $\phi=\psi$.  In this situation, every operator on $\htwo$ intertwines $\tphi$ and $\tpsi$.  

 Suppose $\phi$ is nonconstant and nonunivalent. By the open mapping theorem, there is a nonvoid open set $B$ in $\phiu$ such that $\phi^{-1}(\{b\})$ contains more than one point for each $b\in B$.   Since the set of points in $\U$ at which $\phi'$ vanishes is countable, there  is a $b_0\in B$  such that $\phi^{-1}(\{b_0\})$ contains only points at which $\phi'$ is nonzero.   Let $a_1$ and $a_2$ be distinct points in $\U$ such that $\phi(a_1)=\phi(a_2) = b_0$.  Because $\phi'(a_1)$ and $\phi'(a_2)$ are nonzero, there are disjoint open sets  $A_1$ and $A_2$ in $\U$ with $A_1$  containing $a_1$ and $A_2$ containing $a_2$ such that  $\phi$ restricts to be univalent on each of $A_1$ and $A_2$.  Let $\phi_1$ and $\phi_2$ be the restrictions of $\phi$ to $A_1$ and $A_2$, respectively.    If $\psi$ is such that $\psiu\subset \phi_1(A_1)\cap\phi_2(A_2)$, then $\omega_1:=\phi_1^{-1}\circ \psi$ and $\omega_2:=\phi_2^{-1}\circ\psi$ are distinct holomorphic self-maps of $\U$ such that  $\phi\circ\omega_1 = \phi\circ\omega_2 = \psi$.   It follows that $C_{\omega_1} + C_{\omega_2}$ intertwines $\tphi$ and $\tpsi$.   It is an easy exercise to show that because $\omega_1$ and $\omega_2$ are distinct functions, the intertwining map $C_{\omega_1} + C_{\omega_2}$ is not a weighted composition operator. 
 
  Of course, $C_{\omega_1}$ by itself intertwines $\tphi$ and $\tpsi$, and it is natural to consider the question of whether there are situations where $\tphi$ and $\tpsi$ are intertwined by a linear combination of weighted composition operators but not by a single weighted composition operator.  This cannot happen.  
  
   \begin{theorem}\label{wcos}  The following are equivalent:
   \begin{itemize}
   \item[(a)] $\psi$ is subordinate to $\phi$;
   \item[(b)] $\tphi$ and $\tpsi$ are intertwined by a composition operator;
   \item[(c)]  $\tphi$ and $\tpsi$ are intertwined by a linear combination of weighted composition operators.
  \end{itemize}
  \end{theorem}
\begin{proof} The implications (a)$\implies$(b) and (b)$\implies$(c) hold trivially. We must argue that (c) implies (a).  We prove the contrapositive. Assume $\psi$ is not subordinate to $\phi$; we show inductively there is no linear combination of weighted composition operators intertwining $\tphi$ and $\tpsi$.    Suppose that the single weighted composition operator $C_{\omega, h}$ intertwines $\tphi$ and $\tpsi$; i.e., $C_{\omega, h} \tphi = \tpsi C_{\omega, h}$.  Then for every $f$ in $\htwo$,
   $$
 (\phi\circ \omega)\cdot h \cdot (f\circ\omega) =\psi \cdot h\cdot (f\circ \omega)
   $$
    and it follows (since $h$ is not the zero function) that $\phi\circ\omega = \psi$, contrary to our assumption that $\psi$ is not subordinate to $\phi$.    
  
   Now suppose that, for some positive integer $k$, no linear combination of $k$ weighted composition operators intertwines $\tphi$ and $\tpsi$ yet
  $Y:= \sum_{j=1}^{k+1} C_{\omega_j, h_j}$ intertwines $\tphi$ and $\tpsi$:
 \begin{equation}\label{k1}
 Y\tphi f = \tpsi Yf \ \ {\rm for\ all} \ \ f\in \htwo.
 \end{equation}
 For every nonnegative integer $n$, let $f_n$ be the function defined on $\U$ by $f_n(z) = z^n$. Apply (\ref{k1}) with  $f_n$ replacing $f$, $n =0, 1, \ldots, k$, and do a bit of rearranging to obtain the following system of equations:
\begin{eqnarray*}
\sum_{j=1}^{k+1} h_j \cdot (\phi\circ\omega_j - \psi) &=& 0\\
\sum_{j=1}^{k+1} \omega_j h_j \cdot (\phi\circ\omega_j - \psi) &=& 0\\
\vdots \rule{.5in}{0in}  & = &  \vdots\\
\sum_{j=1}^{k+1} \omega_j^k h_j \cdot (\phi\circ\omega_j-\psi) &=& 0.\
\end{eqnarray*}
For each $z\in \U$, the preceding system may be written in matrix form
\begin{equation}
V(z)u(z) = 0,
\end{equation}
where  $V(z)$ is the Vandermonde matrix generated by $\omega_1(z), \ldots, \omega_{k+1}(z)$, 
$$
V(z) = \left[\begin{array}{cccc} 1 & 1& \cdots  & 1\\ \omega_1(z)&\omega_2(z)&\cdots & \omega_{k+1}(z)\\ (\omega_1(z))^2 &( \omega_2(z))^2 & \cdots & (\omega_{k+1}(z))^2\\ \vdots&\vdots&\cdots&\vdots\\ (\omega_1(z))^k &(\omega_2(z))^k & \cdots & (\omega_{k+1}(z))^k\end{array}\right],
$$
and $u(z)$ is the column vector with entries 
$$
h_j(z)(\phi(\omega_j(z))- \psi(z)) \qquad ( j =1, 2, \ldots k+1).
$$  
If for some $j\in \{1, 2, \ldots, k+1\}$, $h_j$ is the zero function, then $Y$ is a linear combination of $k$ weighted composition operators intertwining $\tphi$ and $\tpsi$, contrary to our inductive hypothesis.   Because no $h_j$ is the zero function, and because $\psi$ is not subordinate to $\phi$, there are uncountably many points $z\in \U$ for which $u(z)$ is not the zero vector.  This means there is an uncountable set  $E$ of points $z\in \U$ at which the Vandermonde matrix $V(z)$ has zero determinant.   Thus, for each $z\in E$  at least two different coordinates in the $(k+1)$-tuple $(\omega_1(z), \omega_2(z), \ldots, \omega_{k+1}(z))$ will be the same. It follows that there exist $j_1, j_2 \in \{1, 2, \ldots, k+1\}$ with $j_1\ne j_2$, such that $\omega_{j_1}(z) = \omega_{j_2}(z)$ for uncountably many $z\in \U$.  This means $\omega_{j_1} = \omega_{j_2}$ and $Y$ is then a linear combination of $k$ weighted composition operators (one of which will be $C_{\omega_{j_1}, h_{j_1} + h_{j_2}}$), contradicting our inductive hypothesis.    Thus, there is no linear combination of $k+1$ weighted composition operators intertwining $\tphi$ and $\tpsi$.  By induction, there is no linear combination of weighted composition operators intertwining $\tphi$ and $\tpsi$, which completes our proof that (c) implies (a).
\end{proof} 

  We conclude this section with two more corollaries of Theorem~\ref{continuation_theorem}.   The first provides additional examples showing that the containment $\psiu\subset\phiu$ is not sufficient for $\tphi\propto\tpsi$.  
 \begin{corollary}\label{p_valentcor} Let $p > 1$.  Suppose that $\phi$ is $p$-valent, that $\phi-\phi(z_0)$ has a zero of order $p$ at some point $z_0\in \U$, and that $\phi$ singly covers some nonvoid open subset $W$  of $\phiu$.  Then  there exists a $\psi$ with $\psiu\subset\phiu$ such that $\tphi\not\propto\tpsi$.
 \end{corollary}
 \begin{proof} Let $w$ be a point in $W$ and let $E$ be a simply connected subregion of $\phiu$ containing both $\phi(z_0)$ and $w$.   Let $\psi$ map $\U$ univalently onto $E$. We claim that $\tphi \not\propto\tpsi$.  Suppose, in order to obtain a contradiction, that $\tphi\propto\tpsi$.  Then by Theorem~\ref{continuation_theorem}, there is a holomorphic self-map $\omega$ of $\U$ such that $\phi\circ\omega = \psi$.   By the choice of $\psi$, there is a point $a\in \U$ such that $\psi(a) = \phi(z_0)$.   Because $\phi$ is $p$-valent and $z_0$ is a zero of $\phi- \phi(z_0)$ of order $p$, $z_0$ is the only point in $\U$ that $\phi$ takes to $\phi(z_0)$.  Thus, using $\phi\circ\omega = \psi$, we see that $\omega(a) = z_0$ and  $\psi'(a) = \phi'(z_0)\omega'(a) = 0$, the latter fact contradicting the univalence of $\psi$. 
 \end{proof} 
  
\begin{figure}[htb]
  \begin{center}
  \includegraphics[height=1.667in, width=1.875in]{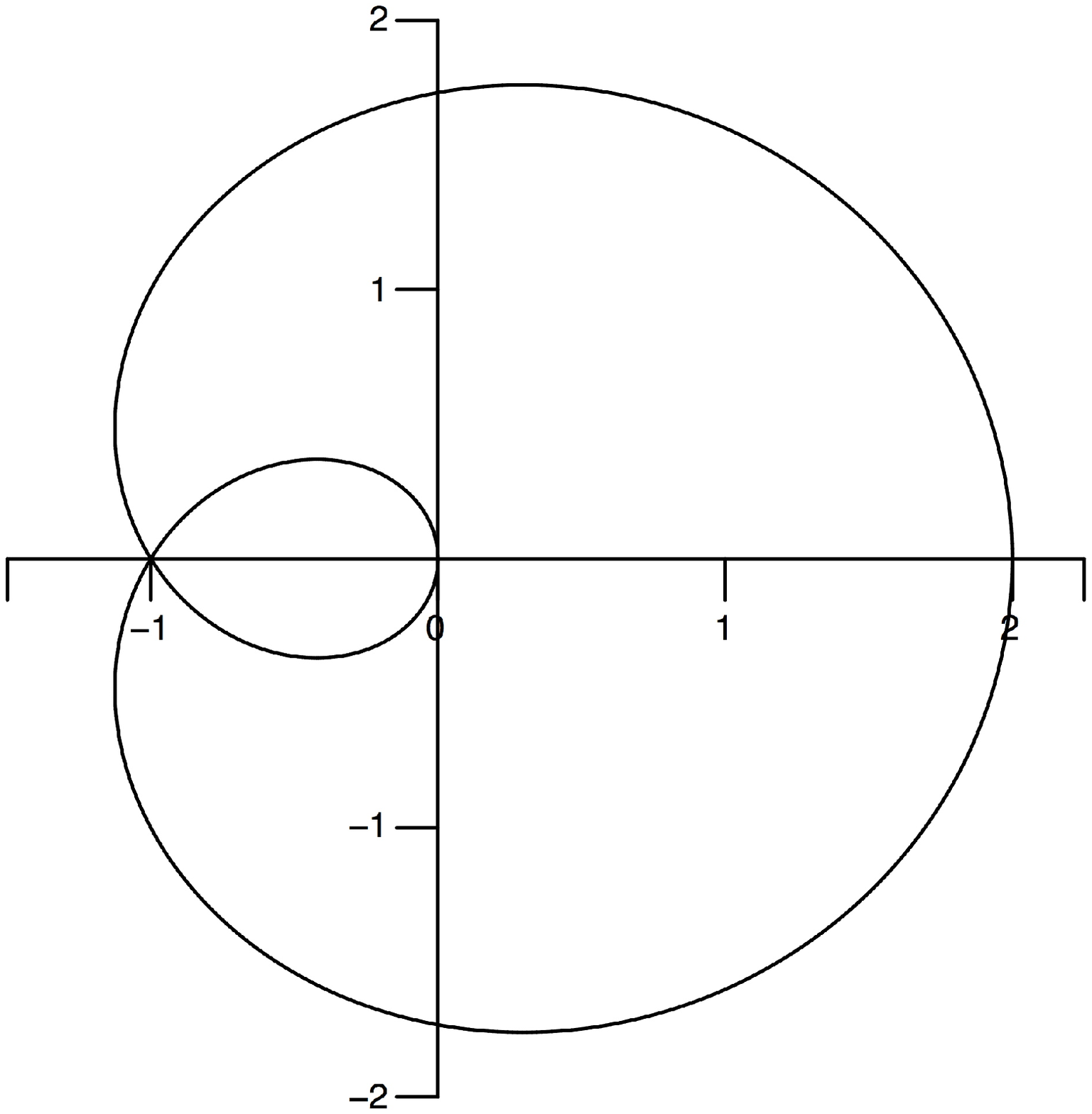}
   \caption{$\phi(\bdu)$ for $\phi(z)=z^2+z$}  \label{figone}
   \end{center}
  \vglue-1em 
\end{figure}

  Note that the proof of the preceding corollary reveals that the hypotheses of $p$-valence and existence of a $p$-th order zero  may be replaced with the single hypothesis that there be a point $w$ in $\phiu$ such that $\phi$ has zero derivative at each $z\in\phi^{-1}(\{w\})$.   For a concrete example  illustrating the corollary, consider the function $\phi(z) = z^2 + z$.  The
mapping properties of $\phi$ are  best understood by writing 
   $$
          \phi(z) =  \left(z+\frac{1}{2}\right)^2-\frac{1}{4}~,
   $$
  whereupon $\phiu$ is revealed to be the region inside the outer cardioid-shaped curve in Figure~\ref{figone}. 
The boundary of $\phiu$  is the image of the larger arc $\Gamma$ of $\bdu$ that lies between the points $e^{\pm i2\pi/3}$, both of which $\phi$ maps to $-1$. The boundary of the small interior loop is the image of the smaller arc $\gamma$ of $\bdu$  between those same two points. The domain $W$ in $\phiu$ that lies outside this smaller loop is singly covered by $\phi$. Indeed,  $V:=\phinv(W)$ is the subdomain  of $\U$ that lies to the right of the reflection of $\gamma$ in the line $\re z = -1/2$; it is taken univalently by $\phi$ onto  $W$.  The points inside the inner loop are doubly covered by $\phi$.  Thus, $\phi$ is at most $2$-valent.   Moreover, $\phi- \phi(-1/2)$ has a zero of order 2. 
 Thus, the hypotheses of Corollary~\ref{p_valentcor} are satisfied, and if we take $\psi$ to be, say, a univalent mapping of $\U$ onto the cardioid $\phi(\U)$, then $\tphi\not\propto\tpsi$, even though $\psiu\subset \phiu$.

 Theorem~\ref{continuation_theorem} also provides a new proof of the following result due to Deddens and Wong \cite[Corollary 4]{DW}.  
 \begin{corollary}  Suppose that $\phi$ singly covers a nonvoid open subset $W$ of $\phiu$; then 
 $$
 \{\tphi\}' = \{T_z\}' .
 $$
 
 \end{corollary}
 \begin{proof}  Let $X$ be a nonzero operator commuting with $\tphi$. Then $X$ intertwines $\tphi$ with itself and we note that the hypotheses of Theorem~\ref{continuation_theorem} are satisfied with $\psi = \phi$.  As we discussed in the paragraph following the proof of Theorem~\ref{continuation_theorem}, $X$ must be a weighted composition operator $C_{\tomega, h}$ where $\tomega$ restricts to $\phinv\circ\psi$ on a nonvoid open set $G$ of $\U$,  but $\phinv\circ\psi$ is the identity function on $G$, so $\tilde{\omega}$ is the identity function on $\U$.  Thus, $X=C_{\tomega, h}$ is the operator of multiplication by $h$ on $\htwo$ so that $h$ must be bounded and  $X = T_h$.  Thus, every operator in the commutant of $\tphi$ is an analytic Toeplitz operator, as desired.
 \end{proof}

\section{Intertwining with $\phi$ Inner} \label{inner_tphi}

The following theorem, originally proved by Deddens \cite[Corollary 2, page 861]{Ded2}, shows that $\tphi\propto\tpsi$ does not imply the containment $\psiu\subset\phiu$: apply the theorem with, e.g.,  $\psi(z)\equiv z$ and $\phi$ any nonconstant inner function  such that $\phiu\neq \U$  (for instance, $\phi$ could be the unit singular  function defined by (\ref{unit_singular}) or more generally the inner function obtained by taking a covering map of $\U$ onto $\U\setminus K$, where $K$ is a compact set of capacity zero \cite[page 37]{CL}). The proof we give below is different from the original one in \cite{Ded2}.

\begin{theorem}\label{intertwining_inner}
Suppose $\phi$ is a nonconstant inner function and $\psi$ is any holomorphic self-map of the unit disc. Then $\tphi \propto \tpsi$. 
\end{theorem}

\begin{proof} Suppose that $\psi$ is constant, taking the value $c\in \U$ at each $z\in \U$.  Because $\phi$ is a nonconstant inner function, Proposition  \ref{pointspec} guarantees that $\sigma_p(\tphistar)$ contains $\U$ (since for any $\alpha\in \U$,  $\phi-\alpha$ has inner factor $(\phi-\alpha)/(1-\bar{\alpha}\phi)$ and outer factor $1-\bar{\alpha}\phi$).  If $X$ is the orthogonal projection of $H^2$ onto the eigenspace for the $\tphistar$ - eigenvalue $\bar{c}$, then it is easy to see that $X$ intertwines $\tphi$ and $\tpsi$, and thus the theorem holds when $\psi$ is constant.

Suppose that $\psi$ is nonconstant.  It is enough to prove the result with $\psi(z)\equiv z$. Suppose for the moment that this has been done, i.e., that we have found a nonzero bounded operator $Y$ on $\htwo$ with $Y\tphi = T_z Y$. For $\psi$ any holomorphic self-map of $\U$, let $X=\cpsi Y$, where $\cpsi$ is the composition operator induced on $\htwo$ by $\psi$. Thus, $X$ is a bounded linear operator on $\htwo$, and since $\psi$ is not constant, $X\neq 0$. Then
   $$
       X\tphi = \cpsi Y \tphi = \cpsi T_z Y = \tpsi\cpsi Y = \tpsi X,
   $$
as desired.

It remains to dispose of the special case $\psi(z)\equiv z$.    Since $\phi$ is a nonconstant inner function,   the range of $\tphi$, which is just $\phi\htwo$, is closed and not all of $\htwo$. Thus, its orthogonal complement, which is the kernel of $\tphistar$, is not the zero-subspace.  Fix a unit vector $f$ in this orthogonal complement and for each non-negative integer $n$ let $f_n=\phi^n f$. One checks easily that our assumptions on $f$ and $\phi$ imply that the sequence $\seq{f_n:n\ge 0}$ is an orthonormal basis for its closed linear span $S$, and that $\tphi f_n=f_{n+1}$ while $\tphistar{f_n} = f_{n-1}$ if $n>0$ and $=0$ if $n=0$. 
Thus, $S$ is a reducing subspace for $\tphi$, with $\tphi$ on $S$ being unitarily equivalent, via the operator $X:S\goesto\htwo$ defined by $Xf_n=z^n~(n=0, 1, 2, \ldots)$, to the forward shift $T_z$ (and $\tphistar$ on $S$ equivalent, via the same operator $X$, to the backward shift $T_z^*$). More precisely, $X\tphi = T_z X$ on $S$, and both $X$ and this operator equation can be extended to all of $\htwo$ by defining $X\equiv 0$ on $S^\perp$.
\end{proof}

\noindent

{\em Remarks.}
(a) Theorem~\ref{intertwining_inner} provides examples of  pairs of Toeplitz operators  that  can be intertwined but cannot be intertwined by composition operators or linear combinations of weighted composition operators.  Given $\phi$ is inner, we need only choose a self-map $\psi$ of $\U$ that is not subordinate to $\phi$.  Then by Theorem~\ref{intertwining_inner}, $\tphi\propto \tpsi$, but by Theorem~\ref{wcos}  no linear combination of weighted composition operators can effect the intertwining.

(b) Note that the intertwinings established by Theorem~\ref{intertwining_inner} are consistent with the necessity part of Deddens' conjecture:
\begin{equation}\label{deddens_containment}
\overline{\psiu}\subset \sigma_p(\tphistar)
\end{equation}
whenever $\phi$ is a nonconstant inner function and $\psi$ is a self-map of $\U$.  In the next section, we show necessity holds when $\psi$ is inner:  if $\tphi\propto \tpsi$ and $\psi$ is inner, then (\ref{deddens_containment}) holds.

(c)  The argument of the first paragraph of the proof Theorem~\ref{intertwining_inner} is a general one: if $\psi$ is a constant function and $\overline{\psiu}\subset \  \sigma_p(\tphistar)$, then $\tphi\propto \tpsi$.  Also applying more generally is the reduction of the second paragraph of the proof, which yields the following transitivity property.

\begin{proposition} \label{subordination_theorem}

   If $\tphi\propto\tpsi$ and $\psi$ is a covering map, then $\tphi\propto T_g$ for every nonconstant holomorphic function $g$ on $\U$ with $g(\U)\subset \psiu$.
\end{proposition}

\begin{proof}
We are assuming that $Y\tphi  = \tpsi Y$ for some bounded linear operator $Y\neq 0$ on $\htwo$. Suppose $g(\U)\subset \psiu$, so that because $\psi$ is a covering map, $g=\psi\circ\omega$ where $\omega$ is a nonconstant holomorphic self-map of $\U$ (see the proof of Theorem \ref{covering_map}). Thus, $X:=\comega Y$ is a bounded operator on $\htwo$ that is not the zero-operator, and
$$
    X\tphi = \comega Y \tphi = \comega \tpsi Y 
              = T_{\psi\circ\omega}\comega Y = T_g X,
$$
i.e., $\tphi\propto T_g$, as desired.
\end{proof}

\section{Intertwining with $\psi$ Inner} \label{inner_tpsi}

We prove that the necessity part of Deddens conjecture is valid, that is $ \overline{\psiu}\subset \sigma_p(\tphistar)$, whenever $\tphi\propto\tpsi$ and $\psi$ is inner.  We obtain this result as a corollary of the following more general theorem.  

\begin{theorem}\label{isometry_thm}
If $S$ is a completely nonunitary isometry on the Hilbert space $\hilb$ and $A\propto S$, then the point spectrum of $A^*$ contains the open unit disc.
\end{theorem}
The definition of ``completely nonunitary,'' as well as the proof of the theorem, depends on the {\em Wold Decomposition} (see, e.g., \cite[Section
23, pp. 111--113]{C}, or \cite[Problem 149, pages 80, 155, and 272]{Hal}).  If $S$ is an isometry on the Hilbert space $\hilb$, then the Wold Decomposition states that $\hilb$ decomposes into the orthogonal direct sum $\calu\oplus\calv$ of two closed subspaces $\calu$ and $\calv$, each invariant for $S$, with $S|_\calu$ unitary and $S|_\calv$ a shift.  More precisely, if $\calw= \ker S^*=(\ran S)^\perp$ then
$$
     \calu = \bigcap_{n=0}^\infty S^n(\hilb) \quad \mbox{and}\quad
      \calv = \infsum \oplus S^n(\calw).
$$
Thus, each vector $v\in\calv$ can be identified with an $\ell^2$ sequence of the form $(w_0,Sw_1,S^2w_2, \, \ldots)$, 
where each $w_j$ belongs to $\calw$, and each such $\ell^2$ sequence corresponds to a vector $v=\infsum S^nw_n$ in $\hilb$. Thus,
the action of $S$ and its adjoint on $\calv$ can be viewed as shifts: thinking of $v=\infsum(S^nw_n)\in\calv$ as a sequence $(w_0, Sw_1, S^2w_2, \, \ldots)$ we have 
$$
   Sv = (0, Sw_0, S^2w_1, S^3w_2, \, \ldots),
$$ 
and because $\sstar S = I$, 
$$   
   \sstar v = ( w_1, Sw_2, S^2w_3, \, \ldots),
$$
where, of course, the ``equality'' in these representations is actually a unitary equivalence. Thus, if $\calu=\seq{0}$ then is $S$  is just a unilateral shift operator or a {\em completely nonunitary  isometry\/}.     
  
It is easy to see (and well known) that a Toeplitz operator induced by a nonconstant inner function is a completely nonunitary isometry.
  
\begin{proposition} Suppose that $\psi$ is a nonconstant inner function; then $\tpsi$ is a completely nonunitary isometry.
\end{proposition}
\begin{proof}  That $\tpsi$ is an isometry on $\htwo$ is clear.  To complete the proof, we show that $\calu:= \cap_{n=0}^\infty \tpsi^n(\htwo) = 
\{0\}$. Suppose $g\in \calu$ and $n$ is an arbitrary positive integer.  Then $g = \tpsi^n w$ for some $w\in\htwo$.  We have  $\|(\tpsistar)^n g\| = \|w\| = \| g\|$.  Since $(\tpsistar)^n$ converges to $0$ strongly (Proposition~\ref{strongconv}), it follows that $g = 0$, completing the proof.
\end{proof} 

Thus, Theorem~\ref{isometry_thm} immediately yields the following.
\begin{corollary}\label{psi_inner} If $\psi$ is a nonconstant inner function and $\tphi\propto\tpsi$, then $\U\subset\sigma_p(\tphistar)$.
\end{corollary}

In addition, to showing that the necessity part of Deddens' conjecture holds when $\psi$ is inner, the preceding corollary yields the following.
\begin{proposition}\label{outer_inner}
 If $\phi$ is outer and $\psi$ is a nonconstant inner function, then $\tphi \not\propto\tpsi$. 
 \end{proposition}
 
 \begin{proof}
 If $\phi$ is an outer function, then the point spectrum of $\tphistar$ does not contain $0$ (Proposition~\ref{pointspec}).  Thus, by Corollary~\ref{psi_inner}, the proposition holds.
  \end{proof}
  
  \noindent
  Let $g$ be a Riemann map of $\U$ onto the right half-disc $\{z\in \U: \re(z) > 0\}$.  Set $\gamma=g^3$ and let $\nu$ be the unit singular function.  Then $\gamma$ is outer and  $\nu$ is inner, and  thus, by Proposition~\ref{outer_inner}, $T_\gamma\not\propto T_\nu$.  This example again shows that $\psiu\subset\phiu$ is not sufficient to yield $\tphi\propto\tpsi$.  Finally, since $T_\nu\propto T_\gamma$ (apply Theorem~\ref{intertwining_inner}), it also illustrates the lack of symmetry of the relation $\propto$ (in a situation where $\phiu = \psiu$). 
    
    We now turn to the proof of Theorem~\ref{isometry_thm}.  This requires a couple of preliminary results. In all that follows, $S$ is a completely nonunitary isometry and we take $\calw=\ker \sstar$, so that we have our Hilbert space $\hilb$ Wold-decomposed as the orthogonal direct sum
$$
     \hilb = \infsum\oplus S^n(\calw)
$$
with $S$ (resp. $S^*$) acting on $\hilb$ as a forward (resp. backward) ``$\calw$-shift,'' as observed earlier.

\begin{definition} 
   For $w\in\calw$ and $\lambda\in\U$, let 
      \begin{equation} \label{kernel_def}
        K(\lambda,w) = \infsum (S^n w) \lambda^n~,
      \end{equation}
   and for $N=0, 1, 2, \, \ldots$\, let
    \begin{equation}\label{kernel_deriv}
       K_N(\lambda,w) = \frac{\partial^N}{\partial\lambda^N} K(\lambda,w)
                               = \infsum c_{N,n} (S^{N+n} w)\lambda^n~,
   \end{equation}
    where
    \begin{equation}\label{cnnf}
   c_{N,n} = (n+1)(n+2)\cdots(n+N).
   \end{equation}  
 \end{definition}

\begin{lemma} \label{kernels}
   For each $w\in\calw$ and $\lambda\in\U$, 
   \begin{equation}\label{eigenvalue}
        \sstar K(\lambda,w) = \lambda K(\lambda,w)
  \end{equation}   
   and more generally, for $N=0, 1, 2, \, \ldots\,$ ,
   \begin{equation}\label{reduce_index}
       (\sstar-\lambda I)K_N(\lambda,w) = NK_{N-1}(\lambda, w)~.
   \end{equation}    
\end{lemma}

\begin{proof}
The key is that since $S$ is an isometry, $\sstar S=I$. Thus, upon applying $S^*$ to both sides of (\ref{kernel_def}), and using the fact that $w\in\ker(\sstar)$ (and noting that the series on the right converges) we easily obtain (\ref{eigenvalue}).

As for  (\ref{reduce_index}), we see from (\ref{cnnf})  that for nonnegative integers $n$ and $N$, \begin{equation}\label{coeff_diff}
     c_{N,n}-c_{N,n-1} = Nc_{N-1,n}.
\end{equation}
Using $S^*S=I$, we obtain
$$
   (\sstar-\lambda I)K_N(\lambda,w) 
          = c_{N,0}(S^{N-1}w) 
             +\sum_{n=1}^\infty(c_{N,n}-c_{N,n-1})(S^{n+N-1}w)\lambda^n
$$
which, in view of (\ref{coeff_diff}) and the fact that
$
     c_{N,0} = N! = Nc_{N-1,0}
$
  yields (\ref{reduce_index}).
\end{proof}

\begin{lemma}\label{density}
 For each $\lambda\in\U$ the collection of vectors $\seq{K_N(\lambda,w): w\in\calw, N=0, 1, 2, \, \ldots\,}$ spans a dense linear subspace of $\hilb$.
\end{lemma}
\begin{proof}
Fix $\lambda_0\in\U$. We wish to show that the set
$$
     \seq{K_N(\lambda_0,w)):N=0, 1, 2, \,\ldots, w\in\calw}
$$
spans a dense linear manifold of $\hilb$. So suppose $h\in\hilb$ is orthogonal to each of the vectors in this set. We wish to show that $h=0$.

To this end, fix $w\in\calw$ and consider the scalar-valued holomorphic function
$F:\U\goesto\C$ defined by
\begin{equation}\label{F_def}
   F(\lambda)=\bilin{K(\lambda,w)}{h} = \infsum\bilin{S^nw}{h}\lambda^n~.
\end{equation} 
Then for $N$ a nonnegative integer, the $N$-th derivative of $F$ is
$$
   F^{(N)}(\lambda) 
      = \bilin{\frac{\partial^N}{\partial \lambda^N}K(\lambda,w)}{h}
      =\bilin{K_N(\lambda,w)}{h}~,
$$
so our hypothesis on $h$ amounts to asserting that
$$
     F^{(N)}(\lambda_0) = 0 \quad\mbox{for}\quad n=0, 1, 2, \, \ldots~,
$$
whereupon, thanks to its analyticity, $F\equiv 0$ on $\U$.  

Thus, each of the power-series coefficients $\bilin{S^nw}{h}$ in (\ref{F_def}) is zero, and since this is true for each $w\in\calw$ we see that $h\perp S^n(\calw)$ for each $n$. Thus, $h\perp\infsum\oplus S^n(\calw)=\hilb$, so $h=0$, as desired.
\end{proof}

  \begin{proof}{\it of Theorem~\ref{isometry_thm}} \quad We are assuming that $XA=SX$ for some linear operator $X\neq 0$. Fix $\lambda\in\U$. 
  
   {\em To show:} $\lambda$ is an eigenvalue for $\astar$. We have, with $\xstar\neq 0$,
\begin{equation} \label{adjointinter}
  \astar\xstar=\xstar\sstar .
\end{equation}
If $\xstar K(\lambda,w)\neq 0$ for some $w\in\calw$, we are done by (\ref{eigenvalue}) of Lemma \ref{kernels}. Otherwise,  by Lemma \ref{density}, there exists a vector $w\in\calw$ and a positive integer $N$, such that 
\begin{equation} \label{nonzero}
    \xstar K_N(\lambda,w) \neq 0.
\end{equation}     
Fix $w$ and choose $N$ to be the {\em least\/} positive integer satisfying (\ref{nonzero}).
Then by (\ref{adjointinter}), $(\astar-\lambda I)\xstar=\xstar(\sstar-\lambda I)$, hence by  (\ref{reduce_index}) of Lemma \ref{kernels} and the minimality of $N$ in (\ref{nonzero}):
$$
    (\astar-\lambda I)\xstar K_N(\lambda,w) = \xstar(\sstar-\lambda I)K_N(w, \lambda)
    =N\xstar K_{N-1}(w,\lambda)=0.
$$
Since $\xstar K_N(\lambda, w)\neq 0$, it is an eigenvector for $\astar$ and $\lambda$ is the corresponding eigenvalue.  
\end{proof}

It's tempting to try to improve  Theorem~\ref{isometry_thm}  by weakening the ``completely nonunitary'' hypothesis to just ``nonunitary.''  This can
not be done, as the following (class of) examples show.

Let $S_0$ be any nonunitary isometry on a Hilbert space $\hilb_0$. Let $U$ be any unitary operator on a Hilbert space $\hilb_1$, and set $S=S_0\oplus U$, a nonunitary isometry on $\hilb_0\oplus\hilb_1$. Let $A=0\oplus U$ and $X=0\oplus I$ on $\hilb_0\oplus\hilb_1$. Then $XA=0\oplus U=SX$, yet the point spectrum of $A^*$ consists of the singleton $\{0\}$ along with the point spectrum of $U^*$ (if there is any), and so lies in the unit circle $\cup\{0\}$.

\section{Extended Eigenvalues}  \label{eevs}
We denote by $\ee{T}$ the collection of {\em extended eigenvalues\/} of the Hilbert-space operator $T$, i.e., the set of complex numbers $\lambda$ such that $\lambda T \propto T$.
Note that since $IT=TI$, we always have $1\in \ee{T}$.    The goal of this section is to examine what our previous results say about the extended eigenvalues of analytic Toeplitz operators $\tphi$, where, as always, $\phi$ is a bounded analytic function on $\U$.  Note that if $\phi$ is a constant function whose value is nonzero then $\ee{\tphi} = \{1\}$, and if the constant value of $\phi$ is zero, then $\ee{\tphi} = \C$.   

{\em For the remainder of this section, we assume $\phi$ is nonconstant.}  Hence, in particular, $\tphi$ is  one-to-one so that $0$ is never an extended eigenvalue.  In fact, 
$$
\ee{\tphi}\subset \{z: |z| \ge 1\} \ \ {\rm for\  every\ (nonconstant)}\ \phi.
$$
The preceding inclusion is a consequence of Proposition~\ref{first_containment}   (Deddens' result labeled D1 in the Introduction), from which it follows that if $\lambda\in \ee{\tphi}$, then $1/\lambda$ multiplies $\phi(\U)$ into its closure. Another consequence of Proposition~\ref{first_containment} is that if $\tphi$ is {\em invertible\/} ($\phi(\U)$ bounded away from $0$), then $\ee{\tphi}$ is a subset of the unit circle, and will just  be the singleton $\{1\}$ if the image of $\phi$ lacks symmetry about the origin.  For example, $\ee{T_{2+z^2}} = \{1\}$.  

When $\phi$ is univalent, Deddens result D2 asserts that necessary and sufficient for $\lambda\in\ee{\tphi}$ is that $1/\lambda$ multiply $\phi(\U)$ to {\em itself}.  Thus, for example, $\ee{T_z} =  \C\setminus \U$, $\ee{T_{z+1}} = [1,\infty)$, and $\ee{T_{z+2}} = \{1\}$ . 

The preceding three examples illustrate how dramatically the collection of extended  eigenvalues for an operator can fail to satisfy a spectral mapping theorem.  On the other hand,  for any operator $T$ and any nonnegative integer $n$, it is easy to see that if $\lambda\in\ee{T}$, then $\lambda^n\in \ee{T^n}$, i.e., $\ee{T}^n\subset\ee{T^n}$.  In fact, Biswas and Petrovic \cite[Theorem 5.2]{BP} have proved  that $\ee{T^n}= \ee{T}^n$.  Note that if $T$ is invertible,  then  $\ee{T^{-1}} = \ee{T}^{-1}$, and hence,  in this case, the equality $\ee{T^n}= \ee{T}^n$ holds for all integers $n$.

Our extensions of Deddens' results yield corresponding extensions of the extended-eigenvalue observations above, namely:
\begin{itemize}
  \item{\em If $\lambda\in\ee{\tphi}$, then for each $a\in\U$, $\phi(a)/\lambda$ is either in $\phiu$ or is capacitarily isolated in its boundary. Thus, if $\partial\phiu$ consists of nondegenerate continua, then $\lambda^{-1}\phiu\subset\phiu$.}
  \item{\em If $\phi$ is a covering map and $\lambda^{-1}\phiu\subset\phiu$, then $\lambda\in\ee{\tphi}$.}
  \item {\em If $\phi$ is a covering map for which $\partial\phiu$ is a union of nondegenerate continua, then} 
  $\lambda\in\ee{\tphi} \iff \lambda^{-1}\phiu\subset\phiu$.
\end{itemize}

Thus, for example, if the image of a covering map $\phi$ is the unit disk with the segment joining $-1/2$ to $1/2$ removed, then $\ee{\tphi} = \{-1,1\}$.  If $\phi$ is a covering map of an annulus centered at the origin having positive inner radius, then $\ee{\tphi}=  \partial \U$.   On the the other hand,  if $\phi$ is  any analytic mapping from $\U$ onto an annulus not centered at the origin, then $\ee{\tphi} =\{1\}$ since, in this situation, $\lambda=1$ is the only scalar for which $\frac{1}{\lambda}\phi(\U)$ is contained in the closure of  $\phi(\U)$.

We conclude with the following corollary of our work in Section 4, which  illustrates the difficulties involved in characterizing $\ee{\tphi}$ when $\phi$ is not, for example, a  covering map or a power thereof.   
\begin{proposition} Suppose that  $\phi(z) = z^2 + z$; then 
$$
\ee{T_{z^2 + z}} = \{1\} \cup \left(\C\setminus\{-4\phiu\} \right).
$$
\end{proposition}
\begin{proof}  The mapping properties of $\phi(z) = z^2 + z$ are discussed in detail following the proof of Corollary~\ref{p_valentcor}. Just as in that discussion,  let $W$ denote the region inside the cardioid $\phiu$ and outside the inner loop of $\phi(\partial \U)$ so that $W$ is singly covered by $\phi$.  Then $\lambda W$ is singly covered by $\lambda \phi$ for any nonzero scalar $\lambda$.  

Suppose that $\lambda\in \ee{\tphi}$; that is,  $T_{\lambda\phi}\propto T_\phi$.  Then $|\lambda| \ge 1$ and it is easy to see that $\phi(\U)$ intersects $\lambda W$ nontrivially.      Thus, we may apply Theorem~\ref{continuation_theorem} to conclude that $\phi$ is subordinate to $\lambda\phi$; that is,  $(\lambda\phi)\circ \omega =\phi$ for some holomorphic self-map $\omega$ of $\U$.  Suppose that there is a number $a\in \U$ such that $\phi(a) = -\lambda/4$.  Then $\phi(\omega(a)) = -1/4$ from which it follows that $\omega(a) = -1/2$.   Hence, $\phi'(a) = \lambda\phi'(-1/2)\omega'(a) = 0$ and we conclude $a = -1/2$ so that $\lambda = 1$.  Thus, if $\lambda\in \ee{\tphi}$, either $\lambda = 1$ or $\lambda \not\in -4\phiu$.  We have shown that $\ee{\tphi} \subset \{1\} \cup (\C\setminus -4\phiu)$.  

To obtain the reverse containment, we work with an extension of $\phi^{-1}$ from $W$ to  $\phiu\setminus (-1,-1/4]$.   Solving $\phi(z) = w$ for $z$, we have 
\begin{equation}\label{two_s}
z = -\frac12 \pm \sqrt{w + \frac14}.
\end{equation}
Take $\sqrt{\cdot}$ to be the principal branch of the square root function, and
 note that for $w\in \phiu\setminus (-1, -1/4]$, $\sqrt{w+1/4}$ must lie in the right half plane.  It follows that $| -1/2 + \sqrt{w + 1/4}\, | <| -1/2 - \sqrt{w + 1/4}\, |$.  Since at least  one of the two solutions given by (\ref{two_s}) must be in $\U$, we conclude that  $\phi^{-1}: = w\mapsto  -\frac12 + \sqrt{w + \frac14}$ takes  $\phiu\setminus (-1, -1/4]$ into $\U$.
 
We know $1\in \ee{\tphi}$. Suppose that $\lambda\not\in -4\phiu$; then $\phi +\lambda/4$ is nonzero on  $\U$.     One sees upon looking at Figure~\ref{figone}  that $\U\subset\phiu\subset2\U$ (or, more rigorously, using the fact that $\phiu$ has polar equation $0\le r<2\cos(\theta/3)$, $|\theta|\le\pi$); hence,  since $0\not \in \phiu + \lambda/4$, we must have $|\lambda| \ge 4$. Thus, $\phi/\lambda$ lies inside $(1/2)\U$;  in particular, it lies in $\phiu$.   Because $\phiu$ is star-like with respect to the origin, the intersection of the nonpositive real axis with any rotate of $\phiu$ is an interval containing $0$.  Since $(1/\lambda)\phiu$ is a rotate of $\phiu$ followed by a ``shrinking'', $(1/\lambda)\phiu$ will also intersect the nonpositive real axis in an interval containing zero.  Since $-1/4$ cannot be in this interval, neither will be any point to its left; thus,  $(1/\lambda)\phiu$ must lie in $\phiu\setminus(-1, -1/4]$.   Letting $\omega = \phi^{-1}\circ (\phi/\lambda)$, we see that $C_\omega$ intertwines $\tphi$ and $T_{\phi/\lambda}$, and hence $\lambda\in \ee{\tphi}$. Thus, any $\lambda\in \C\setminus -4\phiu$ belongs to $\ee{\phi}$, which completes the proof. 
\end{proof}

{\bf Added in proof:}  We recently became aware that results in \cite{BP}, \cite{BPL}, and \cite{K} pertaining to the Volterra operator and its positive integral powers may be extended to its real powers based on earlier work in M. Malamud's papers  \cite{Mal1, Mal2}. For this, as well as clarification of some priority issues, see \cite[page 183 and Remark 2.2 page 194]{Mal3}.


%

\end{document}